\documentclass[article,11pt]{article}
\usepackage[lmargin=25mm, tmargin=25mm, textwidth=160mm, textheight= 221mm]{geometry}
\usepackage{amsmath, amssymb, graphicx, color,bbm}
\usepackage{hyperref}

\usepackage{stmaryrd,amsmath,amssymb, caption,subcaption, geometry,lineno} 

\newcommand{\dsp}{\displaystyle}
\newcommand{\aw}{\operatorname{aw}}
\newcommand{\ecc}{\epsilon} 
\newcommand{\dist}{\operatorname{d}}
\newcommand{\diam}{\operatorname{diam}}
\newcommand{\rad}{\operatorname{rad}}
\newcommand{\N}{\operatorname{N}}

\newtheorem{thm}{Theorem}[section]

\newtheorem{prop}[thm]{Proposition}
\newtheorem{cor}[thm]{Corollary}
\newtheorem{lem}[thm]{Lemma}

\newtheorem{conj}[thm]{Conjecture}

\newtheorem{obs}[thm]{Observation}

\newcommand{\bpf}{ \noindent{\bf Proof:} \newline}
\newcommand{\epf}{\hfill $\square$\newline}

\title{Anti-van der Waerden numbers on Graphs }
\author{Zhanar Berikkyzy\thanks{University of California-Riverside, zhanar@ucr.edu}, Alex Schulte\thanks{Iowa State University, alex.schulte0227@gmail.com}, Elizabeth Sprangel\thanks{Iowa State University, sprangel@iastate.edu},\\ Shanise Walker\thanks{University of Wisconsin-Eau Claire, walkersg@uwec.edu}, Nathan Warnberg\thanks{University of Wisconsin-La Crosse, nwarnberg@uwlax.edu}, Michael Young\thanks{Iowa State University, myoung@iastate.edu}}
\date{\today}

\begin{document}


\maketitle

\section*{Abstract}

\indent \indent In this paper  arithmetic progressions on the integers and the integers modulo $n$ are extended to graphs. This allows  for the definition of the anti-van der Waerden number of a graph.  Much of the focus of this paper is on $3$-term arithmetic progressions for which general bounds are obtained based on the radius and diameter of a graph.  The general bounds are improved for trees and Cartesian products and exact values are determined for some classes of graphs.  Larger $k$-term arithmetic progressions are considered and a connection between the Ramsey number of paths and the anti-van der Waerden number of graphs is established. \\

\noindent {\bf Keywords} anti-van der Waerden number; rainbow; $k$-term arithmetic progression; Ramsey number.

\section{Introduction}

\indent \indent	The motivation for studying the anti-van der Waerden number of graphs originates from extending results on the anti-van der Waerden number of $[n] = \{1,2,\dots,n\}$ and $\mathbb{Z}_n$ to paths and cycles, respectively.  Notice that the set of arithmetic progressions on $[n]$ is isomorphic to the set of non-degenerate arithmetic progressions on $P_n$. Similarly, the set of arithmetic progressions on $\mathbb{Z}_n$ is isomorphic to the set of non-degenerate arithmetic progressions on $C_n$. Therefore, considering the anti-van der Waerden number of $[n]$ or $\mathbb{Z}_n$ is equivalent to studying the anti-van der Waerden number of paths or cycles, respectively.  

The anti-van der Waerden number was first defined in \cite{U}. Many results on arithmetic progressions of $[n]$ and the cyclic groups $\mathbb{Z}_n$ were considered in \cite{DMS}. Other results on colorings (and balanced colorings) of the integers with no rainbow $3$-term arithmetic progressions have also been studied (see \cite{AF}, \cite{AM}, \cite{J}).  In \cite{DMS}, Butler et al.\  determined $\aw(C_n,3)$ as seen in Theorem \ref{cycles}. This theorem was generalized in \cite{finabgroup}. The authors of \cite{DMS} also obtained bounds on $[n]$ and conjectured the result that was proved in \cite{BSY}. This result on $[n]$ is adapted to paths in Theorem \ref{paths}. 


	\begin{thm}\label{paths}\cite{BSY}
If $n \ge 3$ and $7\cdot 3^{m-2} +1 \le n \le 21\cdot 3^{m-2}$, then
	$$\aw([n],3) = \aw(P_n,3) = \left\{\begin{array}{ll} m+2 & \text{ if $n = 3^m$,}\\
											m+3 & \text{ otherwise.}\end{array}\right. $$
									
\end{thm}

In \cite[Theorem 1.6]{DMS} it is shown that $3\leq \aw(\mathbb{Z}_p,3)\leq 4$ for every prime number $p$ and that if $\aw(\mathbb{Z}_p,3)= 4$, then $p\geq 17$. Furthermore, it is shown that the value of $\aw(\mathbb{Z}_n,3)$ is determined by the values of $\aw(\mathbb{Z}_p,3)$ for the prime factors $p$ of $n$.  The notation has been changed, but this result is given in Theorem \ref{cycles}.

\begin{thm}\label{cycles} {\cite{DMS}}
Let $n$ be a positive integer with prime decomposition $n=2^{e_0}p_1^{e_1}p_2^{e_2}\cdots p_s^{e_s}$ for $e_i\geq 0$, $i=0,\ldots,s$, where primes are ordered so that $\aw(\mathbb{Z}_{p_i},3)=3$ for $ 1 \leq i \leq \ell$ and $\aw(\mathbb{Z}_{p_i},3)=4$ for $\ell + 1 \leq i \leq s$. Then, 
		\begin{equation*}
	\aw(\mathbb{Z}_n,3) = 		\aw(C_n,3)=\left\{\begin{array}{ll}
				2 +\sum\limits_{j=1}^\ell e_j + \sum\limits_{j=\ell+1}^s 2e_j & \mbox{if $n$ is odd,} \\
				3 +\sum\limits_{j=1}^\ell e_j + \sum\limits_{j=\ell+1}^s 2e_j & \mbox{if $n$ is even.}
			\end{array}\right.
		\end{equation*}		
	\end{thm}

%
%
%
%
%

%
%
%
%

The anti-van der Waerden number is not a monotone parameter. The inequality $$\aw(C_n,3) = \aw(\mathbb{Z}_n,3) \le \aw([n],3) = \aw(P_n,3),$$ yields examples of this parameter not being monotone.

In this paper, the anti-van der Waerden number of a graph is defined. The definition is inspired by \cite{DMS}, where the anti-van der Waerden numbers were studied over the set of integers $[n]$ and the integers modulo $n$.  First, some fundamental definitions are required.



\indent Let $u,v \in V(G)$, the \emph{distance} from $u$ to $v$, $\dist(u,v)$, is the length of the shortest path from $u$ to $v$. If there is no path from $u$ to $v$, then $\dist(u,v) = \infty$. A \emph{$k$-term arithmetic progression of a graph} $G$, $k$-AP, is a subset of $k$ vertices of $G$ of the form $\{v_1, v_2, \ldots, v_k \}$, where $\dist(v_i,v_{i+1})=d < \infty$ for all $1\leq i<k$.  A $k$-term arithmetic progression is \emph{degenerate} if $v_i=v_j$ for any $i\neq j$.

\indent  An \emph{exact $r$-coloring of a graph} $G$ is a surjective function $c:V(G) \to [r]$.  A set of vertices $S \subseteq V(G)$ is \emph{rainbow} under coloring $c$, if for any $v_i, v_j \in S$, $c(v_i) \neq c(v_j)$ when $v_i \neq v_j$.  Note that degenerate $k$-APs will not be rainbow.  Given a set of vertices $S\subseteq V(G)$, $c(S) = \{c(v) | v \in S\}$, is the set of colors used on the vertices of $S$.


	The \emph{anti-van der Waerden number of a graph} $G$, denoted by $\aw(G,k)$, is the least positive integer $r$ such that every exact $r$-coloring of $G$ contains a rainbow $k$-AP. If $G$ has $n$ vertices and no coloring of $G$ contains non-degenerate $k$-APs, then $\aw(G,k) = n + 1$.

\subsection{Preliminary Results}\label{prelim}

The lower bound for Observation \ref{bounds} is easy to see, and useful, and the upper bound is due to the fact that it is possible to color each vertex of a graph uniquely and fail to have a rainbow $k$-AP.  A \emph{complete graph} on $n$ vertices, written $K_n$, is a graph such that there is an edge between every pair of vertices.  Complete graphs realize the lower bound of Observation \ref{bounds}, namely $\aw(K_n, k) = k$. An \emph{empty graph} on $n$ vertices is a graph on $n$ vertices with no edges.  Empty graphs realize the upper bound of Observation \ref{bounds}. 

\begin{obs}\label{bounds}
Let $G$ be a graph on $n$ vertices, then $k \le \aw(G,k) \le n+1$.
\end{obs}


%
%

A graph $G$ is \emph{connected} if for each pair of vertices $u,v \in V(G)$ there exists a path from $u$ to $v$ in $G$. A graph that is not connected is called \emph{disconnected}. A \emph{connected component} of $G$ is a maximal connected subgraph of $G$. 

\begin{obs}\label{disconnected}
	If $G$ is disconnected with connected components $\dsp\{G_i\}_{i=1}^\ell$, then $$\aw(G,k) = 1 + \dsp\sum_{i=1}^\ell (\aw(G_i,k) - 1).$$
	\end{obs}


Observation \ref{disconnected} is shown by putting $\aw(G_i,k)-1$ colors on each connected component $G_i$ of $G$ and as soon as one additional color is used, one of the connected components contains a rainbow $k$-AP. This observation shows that studying disconnected graphs is equivalent to studying each of its connected components. Therefore, for the remainder of this paper only connected graphs will be considered.\\

In Sections \ref{radsection}, \ref{treesection}, and \ref{graphfamsection} the anti-van der Waerden number on graphs will be considered on $3$-term arithmetic progressions. In Section \ref{radsection},  the radius and diameter of a graph are used to determine bounds on the anti-van der Waerden number of general graphs.  In Section \ref{treesection},  the anti-van der Waerden number on trees is investigated and demonstrates that the anti-van der Waerden number of a tree can be bounded by the number of specific degree two vertices.  In Section \ref{graphfamsection}, the anti-van der Waerden number is determined for some families of graphs and bounds on some graph operations are established. In Section \ref{graphsectionchapter4}, the anti-van der Waerden number on graphs will be considered on $k$-term arithmetic progressions where $k>3$ and connection between the anti-van der Waerden number of a graph and the Ramsey number of multiple paths is established.

\section{Radius and Diameter}\label{radsection}

In this section, anti-van der Waerden numbers are determined and bounded in terms of distance parameters. The \emph{eccentricity} of vertex $v$ in graph $G$ is $\ecc(v) =\max\{\dist(v,u)\,| \, u\in V(G)\}$.  The \emph{radius} of graph $G$ is $\rad(G) =\min\{ \ecc(v) \, | \, v\in V(G)\}$ and the \emph{diameter} is $\diam(G) = \max\{\ecc(v) \, | \, v\in V(G)\}$.  A vertex $v$ is \emph{central} in $G$ if $\ecc(v) = \rad(G)$.

\begin{prop}\label{radprop1}
	If $G$ is a connected graph, then
	
	$$\aw(G,3) \le \left\{\begin{array}{ll}
			\rad(G) + 2 & \text{if $\rad(G) \le 2$},\\
			\rad(G) + 1 & \text{if $3\le \rad(G).$}
	
		\end{array}\right. $$ 
		
\end{prop}
\bpf 
Assume $\rad(G) \le 2$ and let $v$ be a central vertex of $G$. For $1 \le i \le \rad(G)$, let $L_i$ be the set of vertices that are distance $i$ from $v$.  In any exact $(\rad(G)+2)$-coloring of $G$, there exists an $L_i$ that has two colors that are different from the color of $v$.  If $u$ and $w$ are two such vertices in $L_i$, then $\{u, v, w\}$ is a rainbow 3-AP.


%

\begin{figure}[h!]
 	\centering
 	\includegraphics[scale = .4]{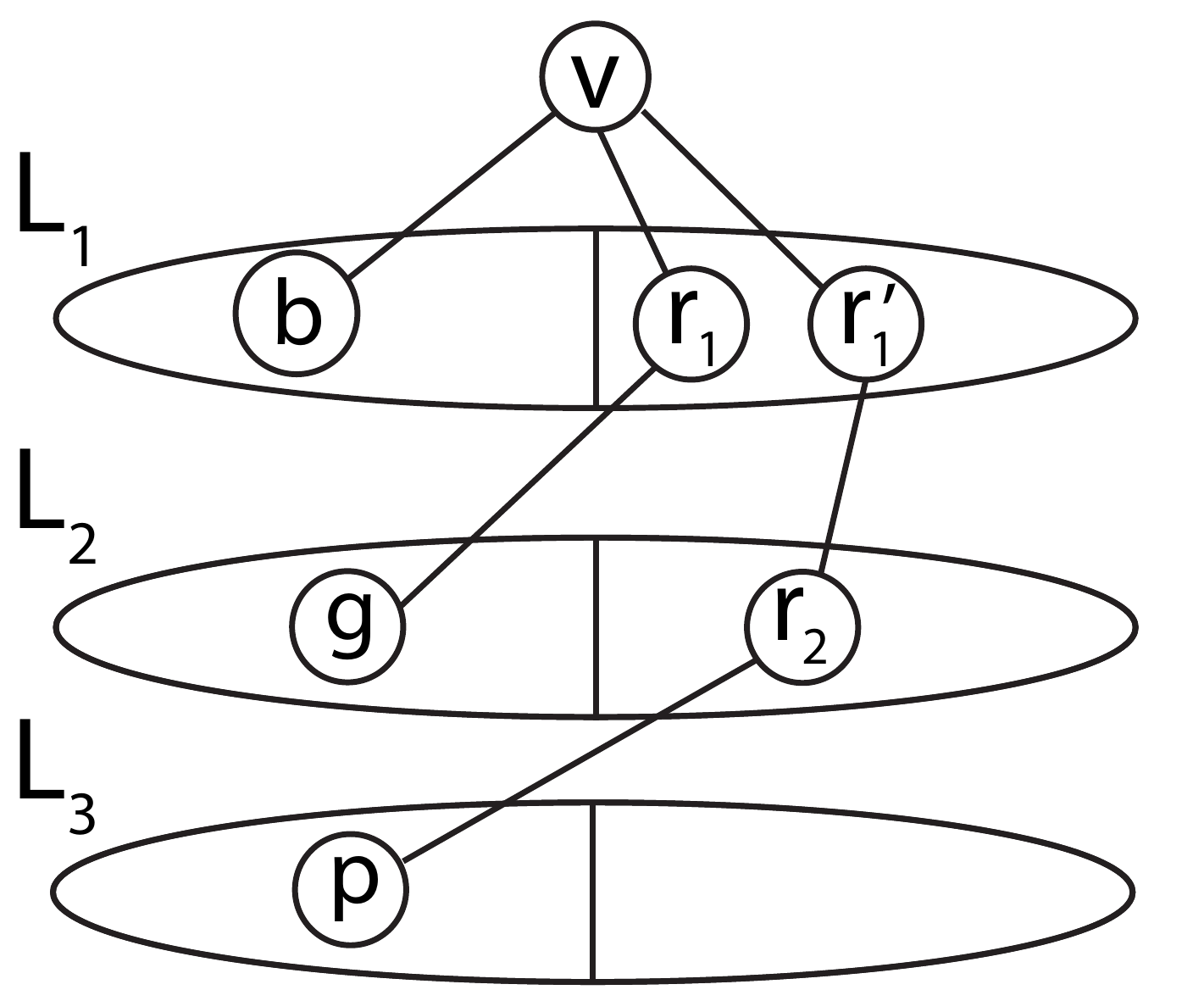}
	\caption{Levels of graph from central vertex $v$, edges drawn must exist.}\label{rad3_1}
 \end{figure}

 Now, assume $G$ has $\rad(G) \ge 3$. Let $c$ be an exact $(\rad(G)+1)$-coloring of $G$ that avoids rainbow $3$-APs.  Without loss of generality, let $c(v) = red$ and $L_i$ be the set of vertices that are distance $i$ from $v$ (see Figure \ref{rad3_1}).  Note that if $v_1,v_2 \in L_i$ where $c(v_1),c(v_2) \neq red$ and $c(v_1) \neq c(v_2)$, then $\{v_1,v,v_2\}$ is a rainbow $3$-AP.  Since there are $\rad(G)$ sets $L_i$, each non-$red$ color appears in exactly one $L_i$.  Let $b\in L_1$, $g \in L_2$ and $p \in L_3$ be vertices that are not colored $red$.  There are no edges between $b$, $g$ and $p$, otherwise there would exist a path on three vertices with three different colors which is a rainbow $3$-AP.  Therefore there exists a vertex colored $red$, $r_1$, in $L_1$ that is adjacent to $g$.  If $r_2$ is adjacent to $b$ or $p$, then there would exist a path on three vertices with three different colors.  This implies there exists $r_2\in L_2$ that is colored $red$ and is adjacent to $p$.  If $r_1$ is adjacent to $r_2$, then $\{b,r_1,p\}$ is a rainbow $3$-AP.  Therefore there exists a vertex colored $red$, $r_1'$, in $L_1$ that is adjacent to $r_2$.  Then either $\{b,r_1',p\}$ or $\{p,b,g\}$ is a rainbow $3$-AP depending on whether or not $br_1'$ is an edge.
 \epf
 
\begin{prop}\label{diam2}
Let $G$ be a connected graph such that $\diam(G) = 2$, then $\aw(G,3) = 3$.
\end{prop}
\bpf  
Since $\diam(G) = 2$, then $\rad(G)$ is either $1$ or $2$.  By Observation \ref{bounds}, $3 \le \aw(G,3)$.  If $\rad(G) = 1$, then $\aw(G,3) = 3$ since there exists vertices $v, u, w$ such that $v$ is a central vertex and $c(v)\neq c(u)$, $c(v)\neq c(w)$ and $c(u)\neq c(w)$, so $\{u,v,w\}$ is a rainbow $3$-AP. Assume $\rad(G) = \diam(G) = 2$ and $c$ is an exact $3$-coloring of $G$.  Let $u$ and $v$ be adjacent vertices such that $c(u) \neq c(v)$.  Let $w$ be a vertex such that $c(w) \neq c(u),c(v)$.  If $w$ is adjacent to $u$ or $v$, there is a rainbow $3$-AP.  If $w$ is not adjacent to $u$ or $v$, then $\{u,w,v\}$ is a rainbow $3$-AP.
\epf


\section{Trees}\label{treesection}

In this section, anti-van der Waerden numbers are bounded for trees.  The fact that connected subgraphs of trees are isometric (distance preserving) is essential for the results in this section.  The relationship between the anti-van der Waerden number of a graph and its isometric subgraph are established in \cite{rsw} (see Lemma 2.1 and Proposition 2.2).

%

%
%
%

\begin{lem}\label{rainbowpath}
	In any coloring of a tree with no rainbow $3$-APs, there exists a path that contains all of the colors.
\end{lem}

\bpf 
Suppose $c$ is an exact $r$-coloring of tree $T$ with no rainbow $3$-APs.  Let $T'$ be the smallest subtree of $T$ that contains all $r$ colors.  If $v$ is a leaf in $T'$, then no other vertex in $T'$ has color $c(v)$; otherwise, $T'$ is not the smallest subtree containing all $r$-colors.  If $T'$ is not a path, then $T'$ has at least three leaves, namely $u$, $v$ and $w$.  Without loss of generality, suppose $\dist(u,v) \le \dist(v,w)$.  Then $\{u,v,x\}$ is a rainbow $3$-AP where $x$ is the vertex on a shortest $vw$-path such that $\dist(u,v) = \dist(v,x)$.  Therefore, $T'$ must be a path.
\epf



%



\begin{prop}\label{pathprop}
	If $T$ is a tree with $\diam(T) =d$, then
	
			$$\aw(T,3) \le \left\{\begin{array}{ll}
						 \aw([d+1],3)  & \text{ if $d\neq 3^m$ for all $m\in \mathbb{Z}$},\\
			                     \aw([d+1],3) + 1&  \text{ if $d = 3^m$ for some $m\in \mathbb{Z}$.}
			                     \end{array}\right.$$

\end{prop}

\bpf

Let $c$ be an exact $(\aw(T,3)-1)$-coloring of $T$ with no rainbow $3$-APs.  By Lemma \ref{rainbowpath}, there exists a path $P$ that contains every color.  Therefore $\aw(T,3)-1 \le \aw([|P|],3) - 1$.  This implies that $\aw(T,3) \le \aw([|P|],3)$.  Since $|P| \le d+1$, then $\aw([|P|],3) \le \aw([d+1],3) + 1$, by Theorem \ref{paths}.  Furthermore, this can be improved to $\aw([|P|],3) \le \aw([d+1],3)$ when $d \neq 3^m$.
\epf

	

%
%

The remainder of this section focuses on the degree two vertices of trees and defines bijacent vertices which can bound the anti-van der Waerden number of a tree.

A \emph{comb} is a graph obtained by adding a leaf to every vertex of a path.  A \emph{broken comb} is a connected subgraph of a comb that has a unique pair of leaves that realize the diameter.  A unique pair of leaves that realize the diameter of a broken comb are called \emph{antipodal vertices}.  A \emph{bijacent vertex} is a degree two vertex with two degree two neighbors.  Observation \ref{closeleaves} gives a reduction technique for computing anti-van der Waerden numbers of graphs with pendant vertices. 

\begin{obs}\label{closeleaves}
	If a graph $G$ has an exact $r$-coloring (with $r \ge 3$) that avoids rainbow $3$-APs and has two leaves that are distance $2$ from each other, then those leaves must of be the same color.  If $u$ and $v$ were leaves distance two from each other and $c(u)\neq c(v)$, then $\{u,x,v\}$ would be a rainbow $3$-AP for any vertex $x$ with $c(x) \neq c(u),c(v)$.  This allows  a graph to be reduced so that no vertex is adjacent to more than one leaf.
\end{obs}

The \emph{neighborhood} of vertex $v$ in graph $G$, denoted $\N_G(v)$ or simply $\N(v)$ when the context is clear, is the set of vertices that are distance $1$ from $v$.

\begin{lem}\label{deg3brokencomb}  
Let $T$ be a broken comb with antipodal vertices $u$ and $v$, and $c$ be an exact $r$-coloring (with $r \ge 3$) with no rainbow $3$-APs such that  $c(u)$ and $c(v)$ are both unique. If  $z\in V(T)$ and $\deg(z)=3$, then $c(z)=c(y)$ for all $y\in N(z)$. 
\end{lem}
\bpf
Consider the path from $u$ to $v$ and label the vertices $u=v_0,v_1,\dots,v_d = v$ where $v_iv_{i+1}\in E(T)$ for $0\le i \le d-1$.  Let $v_i$ be a degree three vertex, and, without loss of generality, assume $i \le d/2$ and  $w \in N(v_i)$ such that $\deg(w) = 1$.  Note that $w$ is not $u$ or $v$, by definition of a broken comb, and so for all $x\in N(v_i)\setminus \{w\}$, either $\{w,u,x\}$ or $\{w,v,x\}$ is a $3$-AP.  Thus, $c(x)=c(y)$ for all $x,y\in N(v_i)$, and so $\{u,v_i,v_{2i}\}$ and $\{u, w, v_{2i}\}$ are $3$-APs.  Therefore, $c(v_i) = c(v_{2i}) = c(w)$ and $v_i$ is the same color as all of its neighbors.
\epf

\begin{lem}\label{bijacentbrokencomb}
Let $T$ be a broken comb with antipodal vertices $u$ and $v$ and $c$ be an exact $r$-coloring (with $r \ge 4$) with no rainbow $3$-APs such that  $c(u)$ and $c(v)$ are both unique. Consider the path from $u$ to $v$ and label the vertices $u=v_0,v_1,\dots,v_d = v$ where $v_iv_{i+1}\in E(T)$ for $0\le i \le d-1$.  
If $c(v_j)\neq c(v_i)$ for all $i< j$ and $j \ge 2$, then $v_j$ is a bijacent vertex. Similarly, if $c(v_j)\neq c(v_i)$ for all $i> j$ and $j \le d-2$, then $v_j$ is a bijacent vertex. 
\end{lem}
\bpf
 Consider the path from $u$ to $v$ and label the vertices $u=v_0,v_1,\dots,v_d = v$ where $v_iv_{i+1}\in E(T)$ for $0\le i \le d-1$.   Let $\alpha \in c(T)\setminus c(\{v_1,u,v\})$ and $v_j$ be the vertex with  the smallest index $j$ such that $c(v_j) =\alpha$.  This implies that $\deg(v_{j-1}),\deg(v_j) = 2$, because if either had degree 3, then $c(v_{j-1}) =c(v_j)$  by Lemma~\ref{deg3brokencomb}, and this contradicts that $j$ was the smallest such index. If $\deg(v_{j+1}) = 3$, then $c(v_j) = c(v_{j+1})$ and either $\{u,v_{j/2},v_j\}$ or $\{u,v_{(j+1)/2},v_{j+1}\}$ is a rainbow $3$-AP since $c(v_{j/2})\neq \alpha$ and $c(v_{(j+1)/2})\neq \alpha$ by the minimality of $j$.  Also, if $\deg(v_{j+1}) = 1$, then $j = d-1$ and either $\{u,v_{j/2},v_j\}$ or $\{u,v_{(j+1)/2},v_{j+1}\}$ is a rainbow $3$-AP.  Therefore, $\deg(v_{j+1}) = 2$ which implies that $v_j$ is a bijacent vertex.  By a similar argument using vertex $v$ instead of vertex $u$, it is shown that if $c(v_j)\neq c(v_i)$ for all $i> j$, then $v_j$ is a bijacent vertex. 
\epf

Theorem \ref{log2ub} uses Lemma \ref{bijacentbrokencomb} by finding isometric broken comb subgraphs within a tree $T$.

\begin{thm}\label{log2ub}
If  $T$ is a tree with $\ell$ bijacent vertices, then 

$$aw(T,3) \leq \left\{\begin{array}{ll}

 4 & \text{if }  \ell=0, \\
  \log_2({\ell})+4 & \text{otherwise}.\\
\end{array} \right.$$
\end{thm}
\bpf
  Let $r=aw(T,3)-1$ and consider an exact $r$-coloring $c$ of $T$ with no rainbow $3$-APs.  If $r\le 2$, the result follows, so assume $r\ge 3$. Let $P$ be the smallest path $v_0, v_1, \dots, v_d$ where $v_i v_{i+1}\in E(T)$ for $0\leq i\leq d-1$ with all the colors of $c(T)$; such a path exists by Lemma \ref{rainbowpath}. Let $T'$ be the broken comb induced by  $P\bigcup_{i=1}^{d-1} N(v_i)$ and suppose that $T'$ has  $\ell'\leq \ell$ bijacent vertices, and antipodal vertices $u=v_0$ and $v=v_d$.  Note that since $P$ is the smallest path containing all the colors of $c(T)$ then $c(u)$ and $c(v)$ are unique in $T'$. 

Suppose $\ell=0$. Assume $r\geq 4$ and let $\alpha\in c(T')\backslash c(\{u, v_1, v\})$. Let $v_j$ be the vertex with the smallest index $j$ such that $c(v_j)=\alpha$. Since $\ell=0$, $v_j$ is not bijacent, which contradicts Lemma~\ref{bijacentbrokencomb}, thus, such an $\alpha$ does not exist and $r=3$. Hence, $aw(T,3)=4$.


For $\ell>0$, proceed by strong induction on  $\ell$. Suppose $\ell=1$ and assume $r\geq 4$. Let $\alpha\in c(T')\backslash c(\{u, v_1, v\})$, and let $v_j$ be the vertex with the smallest index $j$ such that $c(v_j)=\alpha$ and let $v_k$ be the vertex with the largest index $k$ such that $c(v_k)=\alpha$. Note that $j< d/2$, otherwise $\{v_{2j-d},v_j,v\}$ is a rainbow 3-AP. Similarly, $k> d/2$. By Lemma~\ref{bijacentbrokencomb}, both vertices $v_j$ and $v_k$ are bijacent, contradicting $\ell=1$. Therefore, there is no such color $\alpha$, i.e. $r=3$ and so $aw(T,3)=4= \log_2(1)+4$. 

Suppose that $\ell\geq 2$. If $r=3$ then the result follows, so assume $r\geq 4$. Let $j$ be the smallest index such that the set of colors used in $u=v_0, v_1, \dots, v_j$ is $c(T')\backslash c(v)$. Let $k$ be the largest index such that the set of colors used in $v_k, v_{k+1}, \dots, v_d=v$ is $c(T')\backslash c(u)$.  By Lemma~\ref{bijacentbrokencomb}, $v_j$ and $v_k$ are bijacent vertices.  Note that $j< d/2$, otherwise $\{v_{2j-d},v_j,v\}$ is a rainbow 3-AP. Similarly, $k> d/2$. Let $\ell_1$ and $\ell_2$ be the number of bijacent vertices between $u$ and $v_j$ and between  $v_k$ and $v$, respectively. Then, $\ell_1+\ell_2<\ell'\leq \ell$, since $v_j$ and $v_k$ are both bijacent vertices. There are $r-1$ colors that appear on the path between $u$ and $v_j$ because every color except $c(v)$ appears between $u$ and $v_j$. 
Likewise, there are $r-1$ colors that appear on the path between  $v_k$ and $v$. 
If $\ell_1=0$ or $\ell_2=0$, then $r-1=3$ and hence $aw(T,3)=5\le\log_2(\ell)+4$. Suppose that $0<\ell_1$ and $0<\ell_2$. 
By induction, $r-1\leq \log_2(\ell_1)+4$ and $r-1\leq\log_2(\ell_2)+3$ which implies 
$ 2^{r-1-3}\leq \ell_1$ and $ 2^{r-1-3}\leq \ell_2$. Thus, $2^{r-3}\leq \ell_1+\ell_2<\ell'\leq \ell$.  Hence, $r\leq\log_2{\ell'}+3$ and  $aw(T,3)=r+1\leq \log_2{\ell'}+4\leq \log_2{\ell}+4$.
\epf

The authors believe Theorem \ref{log2ub} can be improved.  In particular, there is an upper bound similar to the bound on paths in Theorem \ref{paths}, which is a logarithmic function to base three of the number of bijacent vertices contained in the tree.

\begin{conj}\label{treeconj}
	There exists a constant $C$ such that if $T$ is a tree with at most $\ell$ bijacent vertices in every path, then $\aw(T,3) \le \log_3(\ell) + C$.\\
\end{conj}


\begin{prop}\label{oddtreeprop}
	If $T$ is a tree such that $d= \diam(T)\geq 3$ and $d$ is odd, then $4\le \aw(T,3)$.
\end{prop}

\bpf
	Let $u$ and $v$ be vertices in $T$ such that $\dist(u,v) = d$.  Now an exact $3$-coloring of $T$ that avoids rainbow $3$-APs is described.  First, let $c(u) =  blue$ and $c(v) = green$.  For $w\in V(T)\setminus\{u,v\}$, define
	
	$$c(w) := \left\{ \begin{array}{ll}  green & \text{if } \dist(u,w) = d,\\
							   blue & \text{if } \dist(v,w) = d,\\
							   red & \text{otherwise. } \end{array} \right. $$
							   
	Note that all vertices with eccentricity $d$ are colored $blue$ or $green$.  Furthermore, given a pair of vertices such that one is colored $blue$ and the other is colored $green$, their distance must be $d$.  Since $d$ is odd, there is no rainbow $3$-AP such that the middle vertex is colored $red$.  Therefore, without loss of generality, a rainbow $3$-AP must have the first vertex colored $blue$, the second vertex colored $green$, and the third vertex colored $red$.  However, this would imply a $red$ vertex has eccentricity $d$, which is a contradiction.  Therefore, $4 \le \aw(T,3)$.
\epf

\begin{cor}
If $T$ is a tree such that $\diam(T) \ge 3$ is odd and has no bijacent vertices, then $\aw(T,3)=4$.
\end{cor}

\section{Classes of Graphs}\label{graphfamsection}

In this section, the anti-van der Waerden number is determined for complete binary trees and hypercubes. There are also results on some graph operations.

Let $\mathcal{B}_n$ be the complete binary tree on $2^{n+1}-1$ vertices. Note $\aw(\mathcal{B}_0,3)=2$ since there is only one vertex in $\mathcal{B}_0$ and $\aw(\mathcal{B}_1,3) = 3$ by Proposition \ref{diam2}.

%

       \begin{figure}[h!]
 	\centering
 	\includegraphics[scale = .5]{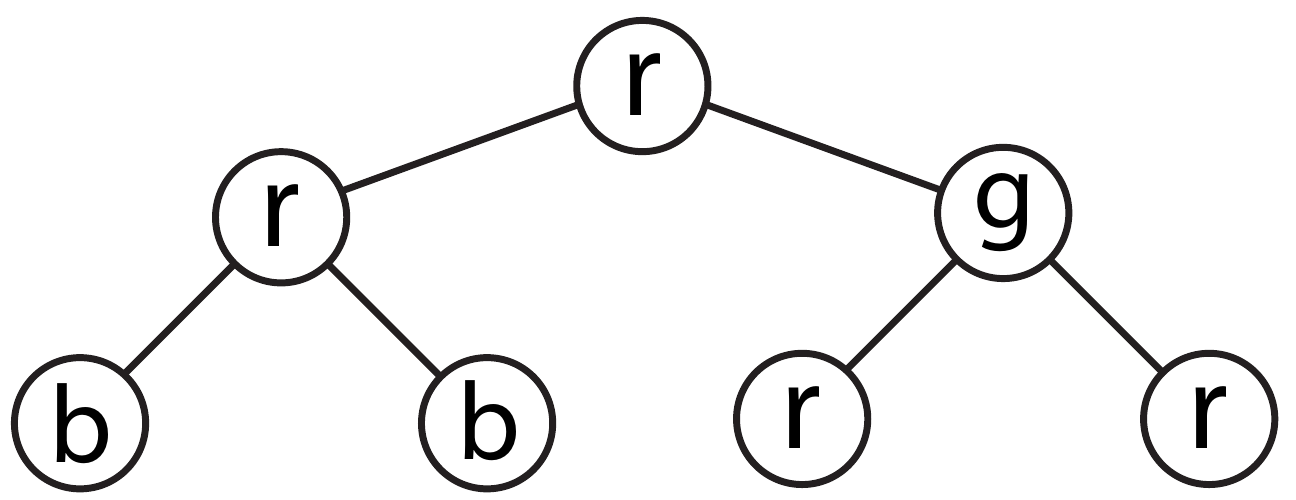}
	\caption{A $3$-coloring of $\mathcal{B}_2$ with no rainbow $3$-AP.}\label{b2fig}
 \end{figure}

\begin{prop}\label{binarytree} For all $n\geq 2$,
$\aw(\mathcal{B}_n,3) = 4$.
\end{prop}
\bpf
For the upper bound, $\aw(\mathcal{B}_n,3) \leq 4$ by Theorem \ref{log2ub} since $\mathcal{B}_n$ has no bijacent vertices.  The case when $n=2$ has a lower bound using the coloring in Figure \ref{b2fig}. Let $n\geq 3$ and $v_0$ be the root of $\mathcal{B}_n$, then following cases give an exact $3$-coloring with no rainbow $3$-APs.  If $n$ is even, let

%

$$c(w) := \left\{ \begin{array}{ll}  red & w = v_0,\\
							   blue & \text{if } \dist(w,v_0) = n-1,\\
							   green & \text{otherwise. } \end{array} \right. $$						   


If $n$ is odd, let

$$c(w) := \left\{ \begin{array}{ll}  red & w = v_0,\\
							   blue & \text{if } \dist(w,v_0) = n,\\
							   green & \text{otherwise. } \end{array} \right. $$

\epf

Let $G$ be a graph with $|V(G)|=n$ and $H$ be a graph with $|V(H)|=m$. The \emph{Cartesian product} of two graphs, denoted $G \square H$, is the graph where $V(G \square H)=V(G)\times V(H)$ and $(u,v)(x,y)\in E(G \square H)$ if and only if $u=x$ and $vy\in E(H)$, or $ux\in E(G)$ and $v=y$. 

It is convenient to consider  $G\square H$ with $m$ copies of $G$, labeled $G_1, G_2, \ldots G_m$, where $G_i$ is the subgraph induced by the vertices of $G\square H$ of the form $(u,i)$ for each $i\in V(H)$.

\begin{lem}\label{diffatmostone}
Let $G$ be a connected graph on $n$ vertices and $H$ be a connected graph on $m$ vertices. Let $c$ be a coloring of $G\square H$ with no rainbow $3$-APs.  If $G_1, G_2, \dots, G_{m}$ are the labeled copies of $G$ in $G\square H$, then $|c(G_j)\setminus c(G_i)|\leq 1$ for all $1\leq i, j \leq m$. 
\end{lem}
\bpf
This proof begins with the special case where $G=P_n$ and $H=P_m$. The  general case is proved by reducing it to the special case. Let $c$ be a coloring of $G\square H$ with no rainbow $3$-APs.\\ 

\textit{Case 1:} $G=P_n$ and $H=P_m$.\\
Assume, for the sake of contradiction, that $|c(G_j)\setminus c(G_i)|>1$. Without loss of generality, suppose $red$ and $blue$ appear in $G_j$ but not $G_i$. Note $i\neq j$, otherwise $|c(G_j)\setminus c(G_i)|=0$. Let $(v,j), (w,j)\in G_j$ such that $c((v,j))=red$ and $c((w,j))=blue$ and $v < w$. If $|j-i|=1$, then $\{(v,j), (w,j), (v',i)\}$ is a rainbow $3$-AP, where $(v',i)$ is the vertex adjacent to $(v,i)$ on the shortest path between $(v,i)$ and $(w,i)$ in $G_i$, a contradiction. Therefore, if $|j-i|=1$, then $|c(G_j)\setminus c(G_i)|\leq 1$.

Now suppose, by way of mathematical induction, that $|c(G_j)\setminus c(G_i)|\leq 1$ whenever $|j-i|<k$. Suppose $|c(G_j)\setminus c(G_i)|>1$ and  $|j-i|= k$.  Define $\dist_G(v,w)=\ell$ and consider the subgraph with corners $(v,i), (w,i), (v,j),$ and $(w,j)$. If $\ell$ were even, then $\{(v,j), (v+ \ell/2, i), (w,j)\}$ would be a rainbow $3$-AP since $c((v+ \ell/2, i))\neq red, blue$, hence $\ell$ is odd. 


Let $(v+1,j)$ and $(w-1,j)$ be the vertices between $(v,j)$ and $(w,j)$ on a shortest path in$G_j$ such that $(v,j)$ is adjacent to $(v+1,j)$ and $(w,j)$ is adjacent to $(w-1,j)$ in $G\square H$. Then, $c((v+1,j-1))$, $c((w-1,j-1))\in \{red, blue\}$ by the respective $3$-APs $\{(v,j), (w,j), (v+1,j-1)\}$ and $\{(w,j), (v,j), (w-1,j-1)\}$. Moreover, $c((v+1,j-1))=red$ by the $3$-AP $\{(v,j), (v,i), (v+1,j-1)\}$ and $c((w-1,j-1))=blue$ by the $3$-AP $\{(w,j), (w,i), (w-1,j-1)\}$ since no vertex in $G_i$ is $red$ or $blue$. Therefore, $\{red, blue\}\subseteq c(V(G_{j-1}))$. However, notice that $|c(G_{j-1})\setminus c(G_i)|\leq 1$ by the induction hypothesis since $|j-1-i|<k$. A contradiction since $red$ and $blue$ are in $c(G_{j-1})$ but not $c(G_i)$. Thus, $|c(G_j)\setminus c(G_i)|\leq 1$ for all $1\le i,j\le m$.\\

%

\textit{Case 2:} Without loss of generality $G\neq P_n$. \\
Let $G_i$ and $G_j$ be two of the labeled copies of $G$ with $1 \le i,j \le m$ and assume $|c(G_j) \setminus c(G_i)|>~1$. Without loss of generality, suppose $red$ and $blue$ appear in $c(G_j)$ but not in $c(G_i)$. Let $c((v,j))=red$ and $c((w,j))=blue$. Let $P$ be the shortest path between $v$ and $w$ in $G$ and let $P'$ be the shortest path between $i$ and $j$ in $H$. 

Consider the isometric subgraph formed by $P\square P'$. Let $P_i$ and $P_j$ be the labeled copies of $P$ from $G_i$ and $G_j$ respectively. Notice $|c(P_j)\setminus c(P_i)| > 1$ and this is again case 1 which implies $P\square P'$ has a rainbow $3$-AP. Notice that since $P$ and $P'$ are shortest paths, distances in the subgraph $P\square P'$ are preserved and therefore any $3$-AP in the subgraph is a $3$-AP in $G\square H$. Since $P\square P'$ has a rainbow $3$-AP, $G\square H$ also has a rainbow $3$-AP, a contradiction.
\epf

Rehm et al.\ showed that the Cartesian product of two graphs has anti-van der Waerden number at most 4 in \cite{rsw}. In their paper they use Lemma \ref{diffatmostone} without including the proof, so the lemma and its proof are included here. Theorem \ref{ghleq4}  was originally a conjecture in this paper, but this paper has been updated to include their result. 


\begin{thm}\label{ghleq4}\cite{rsw}
If $G$ and $H$ are connected graphs and $|G|,|H| \ge 2$, then $\aw(G\square H,3)\leq 4$.
\end{thm}

Another family of graphs investigated are \emph{hypercubes}. The hypercube on $2$ vertices is denoted $Q_1$ and $Q_1 = K_2$. For larger hypercubes, define $Q_n = Q_{n-1} \square K_2$ where $Q_n$ denotes the $n$-dimensional hypercube.  Note that $Q_n$ has $2^{n}$ vertices.  Label the vertices with binary strings of length $n$ such that the distance between two vertices is equal to the number of bits they differ by.  Observe that the diameter of $Q_n$ is $n$.

Let  $\vec{0}_n$ and $\vec{1}_n$ be the vertices of $Q_n$ with all 0s and 1s, respectively. For each vertex $v \in V(Q_n)$, define $|v| := \dist(v,\vec{0}_n)$ and define $\overline{v} \in V(Q_n)$ to be the vertex obtained by switching all $n$ bits of~$v$, i.e. the unique vertex that is distance $n$ from $v$.

Let $L_i = \{ v \in V(Q_n): |v| = i \}$, i.e. the set of vertices with $i$ 1s.

%

\begin{thm}\label{hypercubethm}
For $n \ge 2$,
\[
\aw(Q_n,3) = \left\{
\begin{array}{ll}
3 & \mbox{ if } n \mbox{ is even,}\\
4 &  \mbox{ if } n \mbox{ is odd.} 
\end{array}
\right.
\]
Moreover, if $n$ is odd, then there is a unique exact $3$-coloring of $Q_n$ that avoids rainbow $3$-APs.
\end{thm}
\bpf
This proof is by induction on $n$. It is easy to see that $\aw(Q_2,3) = \aw(C_4,3) = 3$. Assume $3 \le n$, $S_0$ is the subset of vertices in  $V(Q_n)$ with 0 as the first bit, and  $S_1$ is the subset of vertices in  $V(Q_n)$ with 1 as the first bit. Let $L_i$ be the subset of vertices in  $V(Q_n)$ with $i$ bits that are 1s.\\

\textit{Case 1:} $n$ is odd.

By Theorem \ref{ghleq4}, $\aw(Q_n,3) \le 4$, so assume $c$ is an exact $3$-coloring of $Q_n$ with no rainbow $3$-APs. Without loss of generality, either $|c(S_0)| = 2$ and $|c(S_1)| = 1$, or $|c(S_0)| = 2$ and $|c(S_1)| = 2$.

If $c(S_0) = \{red, green\}$ and $c(S_1) = \{blue\}$, then without loss of generality assume $c(\vec{0}_n) = green$. For some $i$ there exists an element $x \in L_i \cap S_0$ such that $c(x) = red$. Then for any $y \in L_i \cap S_1$, the 3-AP $\{x, \vec{0}_n, y \}$ is a rainbow since $c(y) = blue$. This is a contradiction. 

Now assume $c(S_0) = \{red, green\}$ and $c(S_1) = \{red, blue\}$.  There is a vertex $x$ with $c(x) = green$ and $c(\overline{x}) \neq green$. So there is an automorphism of the vertex set such that $c(\vec{0}_n) \neq c(\vec{1}_n)$.

If $c(\vec{0}_n) = green$ and $c(\vec{1}_n) = red$, then a layer $L_i$, with $1 \le i \le n-1$, that contains an element with the color $blue$ must also contain an element with color $red$ in $L_i \cap S_0$. This yields a rainbow 3-AP with $\vec{0}_n$. A similar argument holds if $c(\vec{0}_n) = red$ and $c(\vec{1}_n) = blue$.

If $c(\vec{0}_n) = green$ and $c(\vec{1}_n) = blue$, then for each $i$, $L_i$ is one of the following two types:

\[
\begin{array}{ll}
type \mbox{ } 1: & c(L_i) = \{red\},\\
type \mbox{ } 2: & c(L_i \cap S_0) = \{green\} \mbox{ and } c(L_i \cap S_1) = \{blue\}.
\end{array}
\]

Since $|c(Q_n)|= 3$, a layer of type 1 must exist.  If a layer of type 2 exists, then either there exists an $i$ such that $L_i$ is type 1 and $L_{i+1}$ is type 2, or there exists an $i$ such that $L_i$ is type 2 and $L_{i+1}$ is type 1.  In the first case, if $01\vec{x} \in V(L_{i+1})$, then $\{01\vec{x}, 00\vec{x}, 10\vec{x} \}$ is a rainbow 3-AP since $c(01\vec{x}) = green$, $c(00\vec{x})=red$, and $c(10\vec{x})=blue$. In the second case, if $11\vec{x} \in V(L_{i+1})$, then $\{01\vec{x}, 11\vec{x}, 10\vec{x} \}$ is a rainbow 3-AP since $c(01\vec{x}) = green$, $c(11\vec{x})=red$, and $c(10\vec{x})=blue$.  So there does not exist a layer of type 2 and the only vertices not colored red are $\vec{0}_n$ and $\vec{1}_n$.  This yields a unique coloring of $Q_n$ with exactly 3 colors and no rainbow 3-AP.\\

\textit{Case 2:} $n$ is even.

Assume $Q_n$ is colored with exactly 3 colors and no rainbow 3-AP. If $|c(S_i)| \le 2$, for $0 \le i \le 1$, then the (uniqueness) argument used in \emph{Case $1$} will yield a rainbow $3$-AP for every coloring except the coloring of $Q_n$ with $c(\vec{0}_n)= green$, $c(\vec{1}_n) = blue,$ and the remaining vertices having color $red$. However, since $n$ is even this coloring has a rainbow 3-AP, $\{ \vec{0}_n$, x, $\vec{1}_n \}$, where $x \in L_{n/2}$.

If $|c(S_0)| = 3$, then by induction $S_0$ must have the unique extremal coloring.  In particular, $c(0\vec{0}_{n-1}) = green$, $c(0\vec{1}_{n-1}) = blue$, and the remaining vertices of $S_0$ are colored $red$.


If $c(\vec{1}_n)=blue$, then $\{\vec{0}_n, \vec{0}_{n/2}\vec{1}_{n/2}, \vec{1}_n \}$ is a rainbow $3$-AP.  If $c(\vec{1}_n) = green$, then $\{00\vec{1}_{n-2}$,$01\vec{1}_{n-2}$, $11\vec{1}_{n-2}\}$ is a rainbow $3$-AP.  If $c(\vec{1}_n) = red$ and $x \in L_{n-1}$ such that $c(x) \neq blue$ then, either  $\{0\vec{1}_{n-1}, \vec{0}_n,    x\}$ or  $\{0\vec{1}_{n-1}, \vec{1}_n, x\}$ is a rainbow $3$-AP.  Therefore, $c(L_{n-1}) = \{blue\}$ and either $\{000\vec{0}_{n-3}, 110\vec{0}_{n-3}, 011\vec{0}_{n-3}\}$, $\{11\vec{0}_{n-2}, 10\vec{1}_{n-2}, 00\vec{0}_{n-2}\}$, or  $\{11000\vec{0}_{n-5}, 11011\vec{1}_{n-5}, 00011\vec{0}_{n-5}\}$ is a rainbow $3$-AP when $c(110\vec{0}_{n-3})$ is $blue$, $red$, or $green$, respectively.  Thus any $3$-coloring of $Q_n$ yields a rainbow $3$-AP when $n$ is even. 
\epf

\section{The anti-van der Waerden Number of Graphs with $k\geq 4$}\label{graphsectionchapter4}

Thus far, the majority of the results focus on $k$-APs where $k=3$. In this section, the study of the anti-van der Waerden number of graphs is explored for $k\geq 4$. Recall that determining $3$-term arithmetic progressions of a graph can be simplified by focusing on the central vertex of the progression, then finding equidistant neighbors such that all $3$ vertices have distinct colors. However, in longer arithmetic progressions this technique cannot be used. This section seeks new techniques for these longer progressions.

First observe that $\aw(K_n,k)=k$, since any $k$ vertices with distinct colors form a rainbow $k$-AP (Observation \ref{knk}). However, even an early investigation of the complete bipartite graph $K_{m,n}$ varies from the change between $3$-APs and $4$-APs as seen in Corollary \ref{completebipartitecor}. This example of a variation leads to a more thorough examination of $k$-APs with $k\geq 4$.

\begin{obs}\label{knk}
For $n\geq k\geq 1$, $\aw(K_n,k)=k.$
\end{obs}

Like the investigation of $3$-term arithmetic progressions, the investigation of $k$-term arithmetic progressions begins with an examination of graphs small radii and diameter.

\subsection{Dominating Vertices}\label{domvertsubsection}

 The closed neighborhood of vertex $v$ in $G$ is $\N[v] = \N(v) \cup \{v\}$.  Vertex $v$ in graph $G$ is \emph{dominating} if $\N[v] = V(G)$.

\begin{obs}\label{univert}
	If graph $G$ has a dominating vertex, then $\aw(G,3) = 3$.
\end{obs}


Observation \ref{univert} leads to Conjecture \ref{domconj} when considering $k$-term arithmetic progressions for all $k\geq 3$.

\begin{conj}\label{domconj}
	If $G$ is a graph with a dominating vertex, then $\aw(G,k) \leq k+1.$
\end{conj}

Note that Conjecture \ref{domconj} is true for $k=3$ by Observation \ref{univert} and $k=4$ and $5$,  shown in Proposition \ref{45dom}.

\begin{prop}\label{45dom}
	If $G$ is a graph with a dominating vertex, then $k\le \aw(G,k) \le k+1$ for $k\in\{4,5\}$.
\end{prop}
\bpf
	Let $G$ be a graph with dominating vertex $v$. Note that the lower bounds are a consequence of Observation \ref{bounds}.  First consider $k=4$.  If $|G| \le 4$, then the lower bound is achieved by $K_4$ and the upper bound is achieved, vaccuously by, among others, the star $K_{1,3}$.  Assume $|G| \ge 5$ and let $c$ be an exact $5$-coloring of $G$. Since $v$ is dominating and $G$ has $5$ colors, $v$ must have four neighbors whose colors are distinct from each other and distinct from $c(v)$.  If those neighbors form an independent set, then they form a rainbow $4$-AP. If any pair of them are adjacent, then they form a rainbow $4$-AP that includes $v$. Thus, a rainbow $5$-AP exists in all cases so $\aw(G,4) \leq 5$.

Now, assume $k = 5$.  As in the previous case, the lower and upper bounds can be achieved when $|G| \le 5$.  Assume $|G| \ge 6$ let $c$ be an exact $6$-coloring of $G$.  Then $v$ must have five neighbors, $v_1$, $v_2$, $v_3$, $v_4$, and $v_5$,  whose colors are distinct from each other and distinct from $c(v)$. If $v_1$, $v_2$, $v_3$, $v_4$, and $v_5$ form an independent set, then $\{v_1,v_2,v_3,v_4,v_5\}$ is a rainbow $5$-AP. If exactly one edge exists between the five vertices, say $v_1v_2$, then $\{v_1,v_5,v_4,v_3,v_2\}$ is  a rainbow $5$-AP. 
If there are exactly two edges between the five vertices that share a vertex, say $v_1v_2$ and $v_2v_3$, then $\{v_1,v_4,v_2,v_5,v_3\}$ is a rainbow $5$-AP. Finally, if there are at least two vertex disjoint edges between the five vertices, say $v_1v_2$ and $v_3v_4$, then $\{v_1,v_2,v,v_3,v_4\}$ is a rainbow $5$-AP. All cases yield a rainbow $5$-AP, thus $\aw(G,5) \leq 6$.
\epf

An examination of graphs with dominating vertices would not be complete without also examining stars $K_{1,n}$. In Proposition \ref{starprop}, the anti-van der Waerden number is determined for stars when $k\geq 4$. Note when $k=3$, the anti-van der Waerden number of a star is given by Observation \ref{univert}.

\begin{prop}\label{starprop}
If $4\leq k\le n+1$, then $\aw(K_{1,n},k) = k+1.$
\end{prop}
\bpf 
Let $G = K_{1,n}$ with $4\le k\le n+1$ and $\ell_0$ be the central vertex of $G$.  Label the leaves of $G$ as $\{\ell_1,\ell_2,\dots,\ell_n\}$ and note that $k-1 \le n$.  Define $c: V(G) \to \{0,1,2,\dots,k-1\}$ 
 $$c(\ell_i) = \left\{ \begin{array}{ll} i & \text{if } 0\le i \le k-2,\\ k-1 & \text{otherwise. } \end{array} \right. $$
Then $c$ is an exact $k$-coloring where every $k$-AP has at most $k-1$ colors so $ k+1 \leq \aw(G,k)$.  If $|G| = k$, then Observation \ref{bounds} gives $\aw(G,k) \leq |G| + 1 = k+1$.  Now consider $|G| \geq k+1$ which means $k \leq n$.  Let $c'$ be any exact $(k+1)$-coloring of $G$.  There must exist $k$ colors that appear on the leaves and those leaves form a rainbow $k$-AP, therefore $\aw(G,k) \leq k+1$.  Therefore, $\aw(G,k) = k+1$.
\epf

\subsection{Applications to Ramsey Theory}\label{ramappsection}

\indent In this subsection, connections between the anti-van der Waerden number of a graph and the Ramsey number of multiple paths are identified.  Let $G$ be an edge colored graph and $H$ be a subgraph of $G$, $H$ is \emph{edge monochromatic} if every edge of $H$ is the same color. The \emph{Ramsey number}, $R(k,\ell)$, for positive integers $k$ and $\ell$, is the smallest $n$ such that for every $2$-edge coloring of $K_n$ there exists a monochromatic subgraph isomorphic to $K_k$ in color $1$ or a monochromatic subgraph isomorphic to $K_{\ell}$ in color $2$.  The Ramsey number, $R(k_1, k_2, \ldots, k_r)$, for positive integers $k_i$, is the smallest $n$ such that every $r$-edge coloring of $K_n$ contains an edge monochromatic subgraph isomorphic to $K_{k_i}$ in color $i$ for some $i=1, \ldots, r$.  The Ramsey number, $R(G_1, G_2, \ldots , G_r)$, for graphs $G_i$, is the smallest $n$ such that every $r$-edge coloring of $K_n$ contains an edge monochromatic subgraph isomorphic to $G_{i}$ in color $i$ for some $i=1, \ldots, r$.  Finally, define $R_d(G) = R(G, G, \ldots, G)$ where there are $d$ copies of $G$. 

By forming a complete graph where the color of the edge corresponds to the distance between the two vertices in the original graph, the anti-van der Waerden number can be related to the Ramsey number of paths. This relationship is shown in Theorem \ref{awboundbyramsey}.

\begin{thm}\label{awboundbyramsey}
If $G$ is a graph with $\diam(G)=d$ and $G$ contains at least one $k$-AP, then $k\leq \aw(G,k)\leq R_d(P_k).$
\end{thm}
\bpf
Observation \ref{bounds} gives the lower bound. Suppose $R_d(P_k)=r$ and let $c$ be an exact $r$-coloring of $G$. Define $H$ to be the complete graph $K_r$ formed in the following way.  First, let $v_1, v_2, \ldots, v_r\in V(G)$ such that $c(v_i)\neq c(v_j)$ for all $i\neq j$. Let $v_1, v_2, \ldots, v_r$ be the vertices of $H$. Suppose $c^\prime$ is a coloring of the edges of $H$. Let $c^\prime(v_iv_j)=\dist_G(v_i,v_j)$. Then $|c^\prime|=d$ since $\diam(G)=d$. Any edge monochromatic path of length $k$ in $H$ represents a $k$-term arithmetic progression in $G$. Since each vertex of $H$ is a distinct color under coloring $c$, this edge monochromatic path is equivalent to a rainbow $k$-AP in $G$. Note that such an edge monochromatic path exists since $R_d(P_k)=r$. 
\epf

The following theorem was proved in \cite{gg} and was later cited in \cite{rad}.

\begin{thm}\cite{gg}\label{ramseypathupperbound}
For $n\geq m\geq 2$, $R(P_m,P_n)=n+\left\lfloor\frac{m}{2}\right\rfloor-1.$
\end{thm}

Theorem \ref{awboundbyramsey} gives a general bound for any diameter graph. In the following results, graphs with diameter two are explored. 

\begin{thm}\label{awramconnection}
Let $G$ be a graph such that $\diam(G)=2$, then $k\leq\aw(G,k)\leq k+\left\lfloor\frac{k}{2}\right\rfloor-1.$
Moreover, these bounds are tight.
\end{thm}


\bpf
The bounds are given by Theorem \ref{awboundbyramsey} and Theorem \ref{ramseypathupperbound}.  The rest of the proofs shows that these bounds are tight.

Let $G_1$ be obtained from taking a $K_k$ and removing exactly one edge. Let $c$ be an exact $k$-coloring of $G_1$, then each vertex is a unique color. Clearly, $\diam(G_1)=2$ and there is a rainbow $k$-AP since $G_1$ contains a Hamiltonian path. Therefore, $G_1$ realizes the lower bound  since $\diam(G_1)=2$ and $\aw(G_1,k)=k$.

Define $G_2$ to be the complete bipartite graph $K_{n,n}$ where $n\geq k$ with bipartite vertex sets $X$ and $Y$. Suppose $c$ is a  $\left(k+\lfloor\frac{k}{2}\rfloor-2\right)$-coloring of $G_2$ such that the set of colors used on $X$ is $\{1, 2, \ldots, \lfloor\frac{k}{2}\rfloor-1\}$ and the set of colors used on $Y$ is $\{\lfloor\frac{k}{2}\rfloor, \lfloor\frac{k}{2}\rfloor+1, \ldots, k+\lfloor\frac{k}{2}\rfloor-2\}$.
Notice $c$ avoids rainbow $k$-APs since a rainbow $k$-AP of $K_{n,n}$ must either use only vertices of $X$, only vertices of $Y$, or alternate vertices between $X$ and $Y$. Neither of the first two can form a rainbow $k$-AP since neither $X$ nor $Y$ has $k$ colors. The final $k$-AP can not happen either since $X$ only has $\lfloor\frac{k}{2}\rfloor-1$ colors and it needs to have at least $\lfloor\frac{k}{2}\rfloor$ colors. Therefore, $c$ is an exact $\left(k+\lfloor\frac{k}{2}\rfloor-2\right)$-coloring of $G_2$ which avoids rainbow $k$-APs. Hence, $\aw(G_2,k)\geq k+\lfloor\frac{k}{2}\rfloor-1$. 

Now suppose $G_2$ has $\left(k+\lfloor\frac{k}{2}\rfloor-1\right)$-colors. By Theorem \ref{ramseypathupperbound}, $R(P_k,P_k)=k+\lfloor\frac{k}{2}\rfloor-1$, therefore $\aw(G_2,k)\leq k+\lfloor\frac{k}{2}\rfloor-1$. Therefore, $G_2$ realizes the upper bound  since $\diam(G_2)=2$ and $\aw(G_2,k)=k+\lfloor\frac{k}{2}\rfloor-1$.

\epf

In the proof of Theorem \ref{awramconnection}, it is shown that the complete bipartite graph yields the upper bound, this result is shown in Corollary \ref{completebipartitecor}.

\begin{cor}\label{completebipartitecor}
For $m,n\geq k$, $\aw(K_{m,n},k) = k+\left\lfloor\frac{k}{2}\right\rfloor-1.$
\end{cor}

\section*{Acknowledgements}
 The work of Michael Young is supported in part by the National Science Foundation through grant $1719841$.   

\end{document}